\documentclass[9pt,twocolumn,twoside]{pnas-new}

\templatetype{pnasbriefreport} 

\fancypagestyle{firststyle}{
\fancyfoot[R]{\footerfont \textbf{\today}\hspace{7pt}|\hspace{7pt}\textbf{\thepage}}
\fancyfoot[L]{}
}
\usepackage{amsthm}
\usepackage[capitalise]{cleveref}
\usepackage{overpic}

\DefineSpuriousJournalWord{on}
\DefineSpuriousJournalWord{for}
\DefineSpuriousJournalWord{Its}
\DefineSpuriousJournalWord{its}
\DefineJournalAbbreviation{Annals}{Ann}

\usepackage{tikz}
\usepackage{pgfplots}
\usetikzlibrary{calc}
\pgfplotsset{compat=1.18}
\usepackage{graphicx}
\usepackage{amsmath}
\usepackage{amsfonts}
\usepackage{enumitem}
\usepackage{mathtools}

\usepackage{mathrsfs}  
\graphicspath{{Figure/}}

\newtheorem{theorem}{Theorem}
\newtheorem{lemma}{Lemma}

\newcommand{\R}{\mathbb{R}}

\newcommand{\Fdeg}{\mcal{F}_{\rm{deg}}}

\usepackage{pdfpages}

\usepackage{bm}
\newcommand{\mcal}[1]{{\mathcal{#1}}}

\newcommand{\vect}[1]{{\bm{#1}}}

\title{A zero-one law for one-shot system identification}

\author[a,1]{Nicolas Boull\'e}
\author[b]{Diana Halikias}
\author[c]{Samuel E. Otto}
\author[d]{Alex Townsend}

\affil[a]{Department of Mathematics, Imperial College London, London, SW7 2AZ, UK}
\affil[b]{Courant Institute, New York University, New York, NY 10012, USA}
\affil[c]{Sibley School of Mechanical and Aerospace Engineering, Cornell University, Ithaca, NY 14853, USA}
\affil[d]{Department of Mathematics, Cornell University, Ithaca, NY 14853, USA}

\leadauthor{Boull\'e}

\authorcontributions{Author contributions: N.B., D.H., S.E.O., and A.T. designed research; performed research; analyzed data; and wrote the paper.}
\authordeclaration{The authors declare no competing interest.}
\correspondingauthor{\textsuperscript{1}To whom correspondence should be addressed. E-mail: n.boulle@imperial.ac.uk.}

\keywords{scientific machine learning $|$ system identification $|$ analytic systems $|$ matrix recovery $|$ experimental design}

\begin{abstract}
Can a model be identified from one experiment? We study analytic systems that are linearly parameterized by a combination of prescribed dictionary terms, such as partial differential operators and dynamical systems. For a single input-response pair, recovery is possible exactly when the evaluated dictionary terms are linearly independent. We prove a sharp zero-one law: either no input uniquely determines the coefficients, or almost every random input sampled from a nondegenerate Gaussian measure does. This dichotomy reduces one-shot system identification to a question about degenerate inputs and provides an a posteriori certificate in the noiseless setting. Numerical examples recover dynamical systems, nonlinear partial differential equations, and structured matrix families from single trajectory data, while also detecting when an extra probe is necessary. 
\end{abstract}

\dates{This manuscript was compiled on \today}
\doi{\url{www.pnas.org/cgi/doi/10.1073/pnas.XXXXXXXXXX}}

\begin{document}

\maketitle
\thispagestyle{firststyle}
\ifthenelse{\boolean{shortarticle}}{\ifthenelse{\boolean{singlecolumn}}{\abscontentformatted}{\abscontent}}{}

\dropcap{H}ow many experiments are needed to identify physical laws from data? A fundamental challenge in scientific machine learning~\cite{karniadakis2021physics} is to learn models without requiring an excessive amount of data, which may be costly or difficult to obtain. However, modern scientific machine learning can be strikingly data intensive. For example, recent machine learning weather models are trained on decades of reanalysis data~\cite{price2025probabilistic}, while geoscience models require large geo-datasets~\cite{bergen2019machine}. This makes sample complexity a central issue: how many simulations, trajectories, or experiments are actually needed?

In operator learning, this question has led to theory explaining when solution operators of partial differential equations (PDEs) can be learned efficiently from input-output data~\cite{boulle2023elliptic}. Model identification concerns a different task. Rather than learn a black-box solution operator with a neural network~\cite{li2020fourier}, one seeks a human-interpretable formulation given by the coefficients of a governing law in a prescribed dictionary. While recent years have seen numerous breakthroughs in this area through the introduction of symbolic regression~\cite{bongard2007automated,schmidt2009distilling,udrescu2020ai} and sparsity-promoting algorithms such as SINDy~\cite{brunton2016discovering,champion2019data}, the potential identifiability of governing equations from single trajectory data has remained largely unexplained.

Nearby fields address related but different problems. Operator inference learns reduced surrogate dynamics from simulation data~\cite{kramer2024learning}. Structural identifiability and parameter estimation usually assume the model form and observation protocol are fixed, and then ask which parameters, or functions of parameters, are determined. Algebraic criteria have recently been developed for ordinary differential equations~\cite{hong2020global}. Model recovery also appeared in specialized inverse problems with the establishment of global uniqueness for the Calder\'on problem to recover a conductivity coefficient from boundary measurements~\cite{sylvester1987global}. Our work is motivated by recent developments in scientific machine learning~\cite{brunton2016discovering,champion2019data}, but is genuinely different in its focus. We aim to characterize and quantify the set of inputs that guarantee recovery of the ground-truth model.

\begin{figure}[t]
    \centering
    \includegraphics[width=0.45\textwidth]{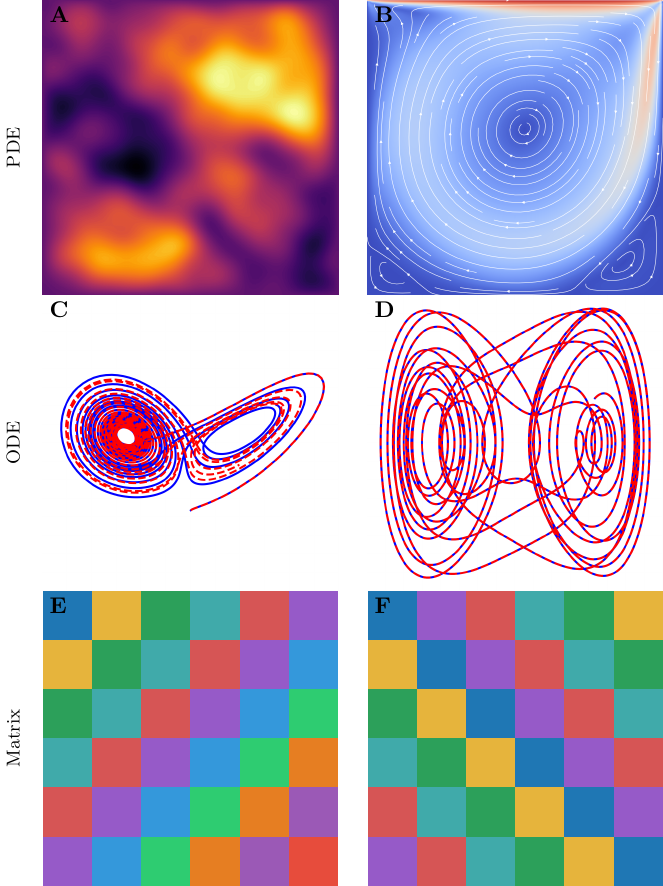}
    \caption{Identifiability across analytic systems. The same rank criterion certifies one-shot recovery for nonlinear differential equations (\textbf{A}, Allen--Cahn; \textbf{B}, Navier--Stokes), dynamical systems (\textbf{C}, Lorenz; \textbf{D}, Duffing), and structured matrix families, while also detecting when an additional probe is needed (\textbf{E}, Hankel) and when one suffices (\textbf{F}, circulant).}
    \vspace{-0.6cm}
    \label{fig:main}
\end{figure}

We consider general analytic systems parameterized by a linear combination of dictionary terms and aim to identify the coefficients from a \emph{single} input-response pair. The unknown object may be a differential equation, a continuous-time dynamical system, or a matrix drawn from a structured linear family (cf.~\cref{fig:main}). These settings all involve a possibly large known dictionary of analytic terms that can be evaluated at a single observed input-response pair, especially for applications where measurement alters the system irreversibly. The problem becomes a coefficient-identification problem where the unknown coefficients must be recovered from that observation, and selecting a good input is crucial to ensure recovery.

\section*{Results}

For one-shot recovery, we assume that one can probe an unknown system with a forcing term $f$ and observe the corresponding response $u$. Our starting point is the observation that coefficient identification reduces to a linear algebra task. We then aim to identify the coefficients in the following analytic system from $(f,u)$:
\begin{equation}\label{eq:model}
    \mcal{L}(u, f) \coloneqq \sum_{i=1}^N c_i^* \mcal{D}_i(u, f) = \mcal{D}_0(u,f),
\end{equation}
where $\mcal{D}_i:\mcal{U}\times\mcal{F}\to\mcal{V}$ are prescribed analytic maps between Banach spaces, $\mathcal{D}_0(u,f)$ is a known term, and the coefficients $c_i^*\in\R$ are unknown. Boundary conditions and initial conditions are built into the spaces $\mcal{U}$ and $\mcal{F}$. Analyticity is a broad structural assumption here. It includes linear maps, polynomial and standard nonlinear dictionary terms, partial derivatives, weak forms of differential operators, and finite-dimensional matrix families. We assume that the solution map $\mcal{G}:\mcal{F}\to\mcal{U}$ associated with \cref{eq:model} is well-posed and analytic and that the input space $\mathcal{F}$ has a countable Schauder basis. This is the infinite-dimensional analogue of a coordinate system, where each input can be expanded in a convergent series in countably many basis functions, and holds for many separable function spaces used in analysis and computation.

Given one input $f\in\mathcal{F}$ to the system modeled by \cref{eq:model} with output $u=\mathcal{G}(f)$, we define an auxiliary parameter-to-recovery map $L_{f,u}:\R^N\to\mathcal{V}$ as $
    L_{f,u}(\vect{a}) =  \sum_{i=1}^N a_i\mcal{D}_i(u,f)$ for $\vect{a}\in\R^N$.
The coefficients $c_i^*$ in \cref{eq:model} are identifiable from the single observation $(f,u)$ if and only if this map is injective. We then aim to characterize the set of degenerate inputs $f\in\mathcal{F}$ for which $L_{f,\mcal G(f)}$ is not injective, defined as
\begin{equation}
    \Fdeg(\vect{c}^*) =
    \left\{ f \in \mcal{F} \mid \exists \vect{a}\in\R^N\setminus\{0\},\, L_{f, \mcal{G}(f)}(\vect{a}) = 0  \right\}.
\end{equation}
To express this question probabilistically, i.e., how likely it is that a randomly chosen input allows for identification, we consider sampling inputs from nondegenerate Gaussian measures on $\mathcal{F}$, so that every nonzero continuous linear observation of the input is Gaussian with positive variance --- a common practice in scientific machine learning~\cite{li2020fourier}. In finite dimensions, this is any Gaussian with nonsingular covariance. The following theorem establishes a zero-one law for one-shot system identification; either model identification is generically possible, or it fails for every input. We phrase this dichotomy in terms of a notion we call compactly slice-null sets, which are sets whose intersections with large enough finite-dimensional affine subspaces translated over any given compact set have zero Lebesgue measure on each slice.

\begin{theorem}[Zero-one law for system identification]\label{thm:analytic_recovery_alternative}
    ~\\
    Exactly one of the following is true:
    \begin{enumerate}[leftmargin=*, label={\roman*)}, nosep]
        \item $\Fdeg(\vect{c}^*)$ is compactly slice-null and $\gamma(\Fdeg(\vect{c}^*)) = 0$ for every nondegenerate Gaussian measure $\gamma$ on $\mathcal{F}$.\label{recov_1}
        \item $\Fdeg(\vect{c}^*) = \mathcal{F}$.\label{recov_2}
    \end{enumerate}
\end{theorem}

Thus, for a fixed analytic system, if a single input succeeds in recovery, the set of degenerate inputs is negligible and recovery is almost surely possible for inputs sampled from a nondegenerate Gaussian measure. On the other hand, if no such input exists and $\Fdeg(\vect{c}^*) = \mathcal{F}$, the obstruction is more fundamental, based on properties like a hidden symmetry or a dictionary relationship that cannot be avoided by changing the forcing. The proof of \cref{thm:analytic_recovery_alternative} studies the prevalence of inputs that make $L_{f,\mcal G(f)}$ degenerate and relies on \cref{lem:level_set}, which extends a standard result in geometric measure theory stating that zero level sets of finite-dimensional analytic functions are negligible~\cite{federer1969geometric}.

\begin{figure}[t]
    \centering
    \includegraphics[width=0.45\textwidth]{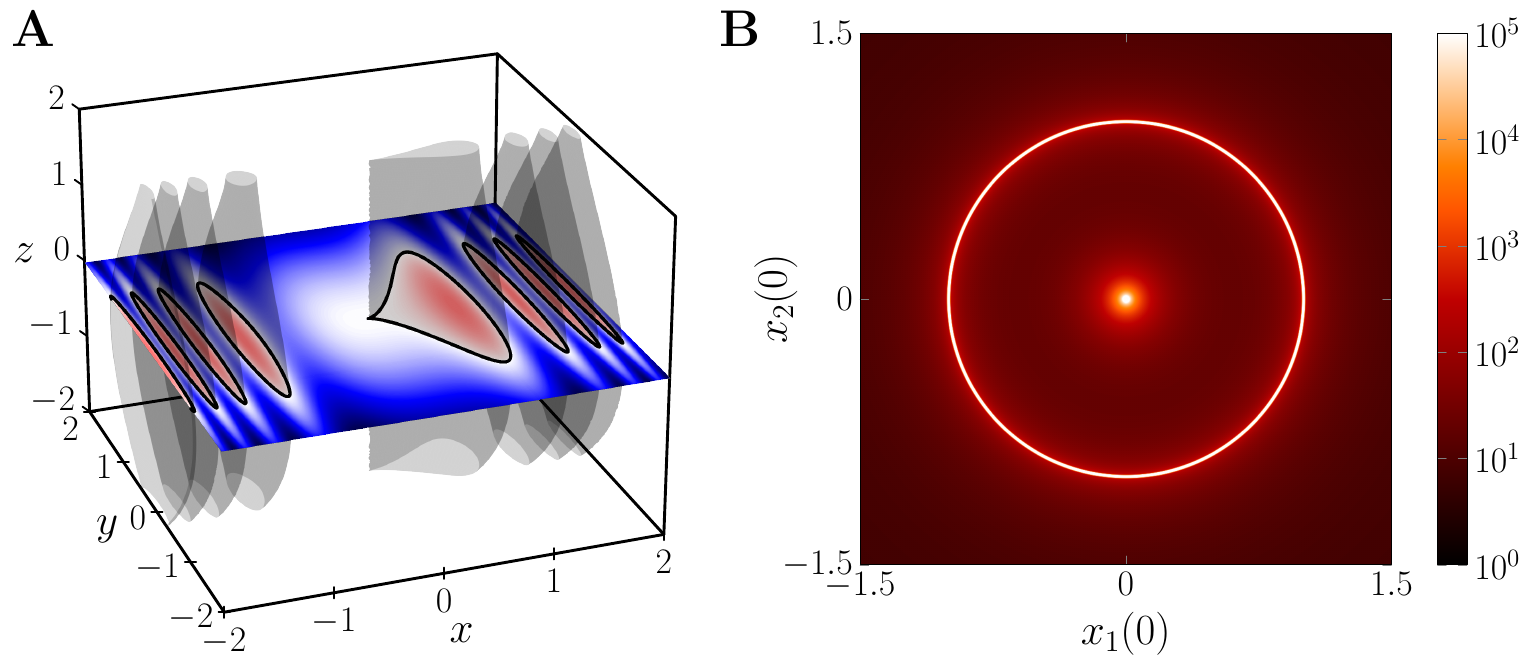}
    \caption{(\textbf{A}) Zero level set of a 3D analytic function illustrated as shaded surfaces. \cref{lem:level_set} ensures that the set is negligible when intersected with a finite-dimensional slice. (\textbf{B}) Condition number of the recovery matrix in \cref{eqn:phi} for a 2D Hopf normal form with a polynomial dictionary as a function of the initial condition.}
    \vspace{-0.6cm}
    \label{fig:level_set}
\end{figure}

\begin{lemma}[Zero level set of analytic functions] \label{lem:level_set}
    ~\\
    Let $\mcal{F}$ be a Banach space with a countable Schauder basis and $\phi: \mcal{F} \to \R$ be a nonzero analytic function. Then, the zero level set of $\phi$ is compactly slice-null and $\phi^{-1}(0)$ has measure zero with respect to any nondegenerate Gaussian measure on Borel sets of $\mcal{F}$.
\end{lemma}

The proof of \cref{lem:level_set} combines a compactness argument with a restriction to finite dimensions to deduce that the zero level set of the analytic function is compactly slice-null, then exploits a disintegration property of Gaussian measures~\cite{bogachev1998gaussian}. To prove~\Cref{thm:analytic_recovery_alternative}, we construct a nonzero analytic function $\phi:\mcal{F}\to\R$ whose zero level set contains  $\Fdeg(\vect{c}^*)$. Suppose that identification succeeds for one input $f_0\in\mathcal{F}$. Then, one can construct linearly independent test functionals $\psi_1,\ldots,\psi_{N}\in\mcal{V}^*$ such that the following analytic function
\begin{equation} \label{eqn:phi}
    \phi:f\mapsto \det\begin{bmatrix}
        \psi_1 (L_{f,\mathcal{G}(f)} e_{1})   & \cdots & \psi_1 (L_{f,\mathcal{G}(f)} e_{N})   \\
        \vdots                                & \ddots & \vdots                                \\
        \psi_{N} (L_{f,\mathcal{G}(f)} e_{1}) & \cdots & \psi_{N} (L_{f,\mathcal{G}(f)} e_{N})
    \end{bmatrix}
\end{equation}
is nonzero at $f_0$ and vanishes on the set of degenerate inputs, where $\{e_1,\ldots, e_N\}$ is the canonical basis of $\R^N$. Applying \cref{lem:level_set} to $\phi$ in \cref{eqn:phi} yields the desired conclusion since the determinant is an analytic function of the input.

We illustrate the zero-one law in \cref{fig:level_set} by plotting the analytic function $\phi$ for a three-dimensional analytic system along with its zero level set and the condition number of the recovery matrix for a 2D Hopf normal form. Here, we verify the conclusions of \cref{lem:level_set} by observing that zero level sets are compactly slice-null. This implies that almost every input enables the identification of the coefficients in these models.

\cref{thm:analytic_recovery_alternative} also provides us a practical identification algorithm that comes with a posteriori certificate for one-shot recovery. After observing the response $u$ to a random Gaussian input $f$, one evaluates the dictionary terms $\mcal{D}_i(u,f)$ through a finite number of test functionals $\psi_1,\ldots\psi_{N}:\mathcal{V}\to\R$ and checks whether the resulting finite-dimensional recovery matrix is full rank, i.e., has nonzero determinant in \cref{eqn:phi}. If this is the case, then the coefficients are uniquely determined by that input, as well as  almost every other random input.

\section*{Discussion}

The zero-one law applies to a wide range of analytic systems, including differential operators, continuous-time dynamical systems, and structured matrix families, as illustrated in \cref{fig:main}.

For differential operators $\mathcal{L}(u) \coloneqq \sum_{i=1}^N c_i^*\mcal D_i(u) = f$, where the dictionary terms only depend on the state variable and partial derivatives of $u$, one can construct a smooth solution $u_0$ such that the evaluated dictionary terms $\mcal{D}_i(u_0)$ are linearly independent. Selecting $f_0 = \mathcal{L}(u_0)$ yields a nondegenerate input that identifies the coefficients. Analytic well-posedness then puts the system in alternative \ref{recov_1} of \cref{thm:analytic_recovery_alternative}: almost every random forcing identifies the coefficients. This covers weak and strong formulations and is illustrated by the Allen--Cahn example in \cref{fig:main}A.

For analytic continuous-time dynamical systems $\dot{u}(t)=\sum_{i = 1}^{N} c^*_i D_i(t, u(t), f(t))$, the input may be an initial condition, a time-dependent forcing, or another experimentally chosen parameter. The zero-one alternative can go either way, but once one trajectory is found to identify the coefficients, almost every random input does. Thus even the chaotic Lorenz system (\cref{fig:main}C) can be identified from a single trajectory in our dictionary model; the Duffing example (\cref{fig:main}D) similarly recovers the governing second-order equation from one observed trajectory.

The same principle also extends to matrix recovery. Suppose an unknown matrix belongs to the family $A(c^*)=\sum_{i=1}^N c_i^*A_i$ for matrices $A_i\in \R^{m\times n}$, and is observed only through matrix-vector products $A(c^*)x$ and $A(c^*)^\top y$. For any fixed collection of probes, identifiability is based on the injectivity of a finite-dimensional recovery map. Circulant matrices fall in the recovering alternative with one random probe, while Hankel matrices do not; stacking two probe-response pairs gives a new block recovery map that is full rank.

\cref{thm:analytic_recovery_alternative} also applies to multi-shot recovery, where several input-response pairs are observed. In this case, the recovery map is a block map that stacks the individual maps for each input-response pair. The zero-one law then applies to this concatenated map, and the same dichotomy holds: either no collection of inputs identifies the coefficients, or almost every random collection does. Thus, we can study multi-shot recovery with the same theory by using the zero-one law applied to a concatenated observation space and increase the number of observations until the recovery map is full rank.

The set of degenerate inputs can be precisely analyzed in specific cases. For single-variable polynomial ODEs, where the dictionary consists of polynomials of $u$ and its derivatives, we find that the degenerate set $\Fdeg(\vect{c}^*)$ does not contain rough functions that are non-differentiable at a dense set of points. On the other hand, the degenerate set is trivial for uniformly elliptic PDEs with known higher-order terms.

In summary, the zero-one law provides a sharp characterization of one-shot system identification for analytic systems. It shows that the success of model recovery from a single experiment is determined by the existence of a single informative input, and that almost every random input will also be informative. This result has practical implications for experimental design and model discovery, as it shifts the focus from collecting many similar snapshots to designing experiments that break dictionary symmetries and using the observed linear system rank as a certificate of success.

Our analysis is idealized in that it assumes exact model evaluation and noiseless observations. In practical settings with discretization and measurement error, near-degeneracy can still be detected through the smallest singular value of the recovery matrix, and multiple experiments can be combined to improve conditioning. Extending the zero-one framework to quantitative stability bounds under noise and to adaptive experiment design is a natural next step. Other open questions include the selection of test functionals and recovery from partial observations of the system.

\matmethods{Proofs of the zero-one law, the compactly slice-null characterization of degenerate sets, and the specialized matrix, dynamical-system, and differential-operator results are provided in SI Appendix. Computationally, one evaluates the dictionary on the observed experiment, tests the resulting finite-dimensional recovery matrix for full rank, and solves the corresponding linear system after rank-revealing column-pivoted QR reduction. This rank check is an a posteriori certificate of uniqueness; with noisy or discretized data, the smallest singular value gives a stability diagnostic. Numerical examples use synthetic input-response data generated from the stated model classes and assemble either weak finite-element or temporal test-function systems.}

\matmethods{
    A key technical result behind the zero-one law is \cref{lem:level_set}, whose proof consists of applying the characterization of root sets of analytic functions~\cite{federer1969geometric} to a finite dimensional auxiliary map. To obtain the second part, we rely on a disintegration property of Gaussian measures on Banach spaces~\cite{bogachev1998gaussian}.
    For each experiment in \cref{fig:main}, we observe one input-response pair $(f,u)$, evaluate the prescribed dictionary terms, and assemble the finite-dimensional recovery matrix by applying selected test functionals to these evaluations. Coefficients are estimated by solving the resulting linear system using rank-revealing column-pivoted QR factorization (see SI Appendix).}

\showmatmethods{} 

\dataavail{Details of the numerical experiments are available in SI Appendix. Code and datasets are publicly available on GitHub at \url{https://github.com/NBoulle/wias}. Proof of the zero-one law and additional results are available in the companion manuscript~\cite{boulle2026one}.}

\acknow{The work of A.T. was supported by the Defense Advanced Research Projects Agency (DARPA) through The Right Space (TRS) Disruption Opportunity (DARPA-PA-24-04-07) and National Science Foundation CAREER grant DMS-2045646.}

\showacknow{} 

\bibliography{pnas-sample}



\section*{Supplementary material}

## Supporting Information Text

## Contents



## 1. Methods

The main contribution of this article is Theorem 1 from the main text which provides a zero-one law for the identifiability of linearly parameterized operators. In this supplementary text, we explain how this result leads to a weak identification procedure for linearly parameterized analytic operators, including settings in which the residuals are tested through functionals rather than inner products, and can be viewed as a generalization of Weak SINDy (1, 2) to arbitrary Banach spaces. We begin with a Hilbert-space formulation leading to a simple and computable recovery algorithm and consider an operator equation of the form:

$$\mathcal{L}(u,f)=\sum_{i=1}^{N} c_i^*\mathcal{D}_i(u,f)=\mathcal{D}_0(u,f), \qquad [1]$$

where $\boldsymbol{c}^*=(c_1^*,\ldots,c_N^*)\in\mathbb{R}^N$ is an unknown coefficient vector, $u\in\mathcal{U}$ and $f\in\mathcal{F}$ belong to Banach spaces. Here, $\mathcal{D}=\{\mathcal{D}_i\}_{i=0}^N$ is a prescribed dictionary of analytic operators $\mathcal{D}_i:\mathcal{U}\times\mathcal{F}\to\mathcal{V}$, where $\mathcal{V}$ is a residual space. Given an observed input-output pair $(f,u)$, we aim to determine whether the functions $\mathcal{D}_1(u,f),\ldots,\mathcal{D}_N(u,f)\in\mathcal{V}$ are linearly independent. When $\mathcal{V}$ is a real Hilbert space, this can be diagnosed through the associated Gram matrix defined as $G_{i,j}=\langle\mathcal{D}_j(u,f),\mathcal{D}_i(u,f)\rangle_{\mathcal{V}}$ for $1\le i,j\le N$. When $G$ is invertible, the dictionary terms are linearly independent and the coefficient vector $\boldsymbol{c}^*$ is uniquely determined from the single pair $(f,u)$ by one linear solve against the right-hand side $b=[\langle\mathcal{D}_0(u,f),\mathcal{D}_1(u,f)\rangle_{\mathcal{V}}\ldots\langle\mathcal{D}_0(u,f),\mathcal{D}_N(u,f)\rangle_{\mathcal{V}}]^\top$ as

$$\boldsymbol{c}^*=\begin{bmatrix}\langle\mathcal{D}_1(u,f),\mathcal{D}_1(u,f)\rangle_{\mathcal{V}} & \cdots & \langle\mathcal{D}_N(u,f),\mathcal{D}_1(u,f)\rangle_{\mathcal{V}}\\ \vdots & \ddots & \vdots\\ \langle\mathcal{D}_1(u,f),\mathcal{D}_N(u,f)\rangle_{\mathcal{V}} & \cdots & \langle\mathcal{D}_N(u,f),\mathcal{D}_N(u,f)\rangle_{\mathcal{V}}\end{bmatrix}^{-1}\begin{bmatrix}\langle\mathcal{D}_0(u,f),\mathcal{D}_1(u,f)\rangle_{\mathcal{V}}\\ \vdots\\ \langle\mathcal{D}_0(u,f),\mathcal{D}_N(u,f)\rangle_{\mathcal{V}}\end{bmatrix}=G^{-1}b.$$

If the Gram matrix is rank deficient, then the observed pair does not uniquely identify the coefficient vector; additional data or an additional normalization convention is required.

This approach also extends to settings in which the residuals live in a general Banach space $\mathcal{V}$, where there is no inner product. We choose linearly independent test functionals $\psi_1,\ldots,\psi_{r_0}\in\mathcal{V}^*$ with $r_0\ge N$ and probe the dictionary terms through the pairings $\psi_i(\mathcal{D}_j(u,f))$. We define the measurement map $\Psi:\mathcal{V}\to\mathbb{R}^{r_0}$ as

$$\Psi(g)=\big(\psi_1(g),\ldots,\psi_{r_0}(g)\big),\quad g\in\mathcal{V}.$$

For the observed pair $(f,u)$, the tested dictionary matrix is $G_{i,j}=\psi_i(\mathcal{D}_j(u,f))$ for $1\le i\le r_0$, and $1\le j\le N$. If $G$ has full column rank, then the tested dictionary terms are linearly independent and the coefficients are uniquely determined from these tests. Conversely, full column rank is equivalent to independence in the residual space provided that $\Psi$ is injective on $\operatorname{span}\{\mathcal{D}_1(u,f),\ldots,\mathcal{D}_N(u,f)\}$.

In practice, we employ a column-pivoted QR factorization $GP=QR$ for numerical stability to detect the numerical rank $r$ of the tested dictionary matrix $G$, remove linearly dependent columns, and solve the resulting reduced triangular system by back substitution. This is important as the dictionary may contain redundant terms. If $G$ is rank-deficient, the original coefficients are not uniquely identifiable as multiple coefficient vectors yield the same operator. In this case, the algorithm returns one representative coefficient vector for the reduced dictionary. Note that the formulation in Eq. (1) may also include spurious dictionary terms that are not present in the true operator $\mathcal{L}$, and for which the associated coefficients are zero. The procedure is summarized in Algorithm 1. Unlike SINDy and Weak SINDy algorithms (1–4), we do not require any sparsity-promoting regularization to recover the coefficients, as we assume that we are in a noise-free regime to verify the applicability of Theorem 1. In the next section, we employ the weak identification algorithm to recover the coefficients of several analytic systems from a single input-output pair. In these examples, the weak formulation improves numerical stability by avoiding direct numerical differentiation of the observed data, as also observed in (1, 5).

**Algorithm 1** Weak Identification of Analytic Systems

**Require:** Source term $f$, observed solution $u$, dictionary $\{\mathcal{D}_i\}_{i=1}^N$, dual space family $\{\psi_i\}_{i=1}^{r_0} \subset \mathcal{V}^*$

Compute the $r_0 \times N$ matrix $G_{i,j} = \psi_i(\mathcal{D}_j(u, f))$

Compute the right-hand side $b = [\psi_1(\mathcal{D}_0(u, f)) \ \ldots \ \psi_{r_0}(\mathcal{D}_0(u, f))]^\top$

Compute a column-pivoted QR factorization $GP = QR$

Determine the numerical rank $r$ from the diagonal of $R$

Keep the independent columns indexed by the first $r$ pivots and remove the dependent columns

Form the reduced upper-triangular system $R_{1:r,1:r} y = Q_{:,1:r}^\top b$

Recover the coefficient vector by undoing the permutation: $c = P \begin{bmatrix} y \\ 0 \end{bmatrix}$

**return** $c$

## 2. Numerical examples

This section illustrates the recovery of analytic systems from a single input-output pair using Algorithm 1 on the six examples presented in Figure 1 of the main text.

**A. Allen–Cahn equation.** We first consider the steady-state cubic-quintic Allen–Cahn equation modeling phase separation (6–9), which is given in strong form by

$$\frac{1}{4}\Delta u - 5u - u^3 + u^5 = f, \qquad u|_{\partial\Omega} = 0, \qquad [2]$$

where $\Omega = [0, 1]^2$, $u \in \mathcal{U} = H^2(\Omega) \cap H_0^1(\Omega)$, and $f \in \mathcal{F} = L^2(\Omega)$. We discretize Eq. (2) using piecewise linear continuous Lagrangian finite elements in the Firedrake finite element software (10) on a $100 \times 100$ mesh. The weak formulation is obtained by multiplying Eq. (2) with a test function $\phi \in \mathcal{V} = H_0^1(\Omega)$ and integrating by parts, resulting in the following weak formulation:

$$-\frac{1}{4}\langle \nabla u, \nabla \phi \rangle - 5\langle u, \phi \rangle - \langle u^3, \phi \rangle + \langle u^5, \phi \rangle = \langle f, \phi \rangle.$$

where the boundary term vanishes because the test function has zero trace on the boundary of $\Omega$. We solve the resulting equation using a random source term $f \sim \mathcal{GP}(0, K)$, with squared-exponential covariance kernel $K(x, y) = \sigma^2 \exp(-\|x - y\|^2/2\ell^2)$, with variance $\sigma^2 = 10$ and length scale $\ell = 0.05$. To sample such a forcing term in the finite-dimensional discretization, we evaluate the kernel at the finite element nodes, form the covariance matrix $K_h$, compute a Cholesky factorization $K_h = LL^\top$, and set $f_h = L\xi$ with $\xi \sim \mathcal{N}(0, I)$. The nonlinear system is then solved with Newton's method using backtracking line search with absolute tolerance of $10^{-12}$, while the linearized system is solved using a sparse LU direct solver (11).

To apply Algorithm 1, we choose test functionals from the finite-element discretization as follows: let $\{\varphi_i\}_{i=1}^{r_0}$ be the nodal basis of the finite element space, and define $\psi_i : w \mapsto \langle w, \varphi_i \rangle$. On a $100 \times 100$ element mesh there are 101 nodes per coordinate direction, hence $r_0 = 101^2$ nodal basis functions. We use the dictionary $\{u_{xx}, u_{xy}, u_{yx}, u_{yy}, 1, u, u^2, u^3, u^4, u^5\}$ consisting of second-order derivatives and monomials up to degree five with right-hand side $\mathcal{D}_0(u, f) = f$. Note that this dictionary contains both spurious terms and the redundant mixed derivative terms $u_{xy} = u_{yx}$, which are removed by the QR factorization in Algorithm 1. We then assemble the weak forms against all basis functions, leading to a system of $r_0 = 101^2$ equations for the $N = 10$ coefficients in the dictionary, and solve the system using a column-pivoted QR factorization. The recovered coefficients match the ground truth ones with maximum absolute error of $1.5 \times 10^{-12}$.

**B. Navier–Stokes equations.** In the second example, we consider the incompressible Navier–Stokes equations for the steady lid-driven cavity flow on the unit square $\Omega = [0, 1]^2$:

$$\begin{aligned} -\nabla \cdot (2\nu\epsilon(u)) + (u \cdot \nabla)u + \nabla p &= 0, \\ \nabla \cdot u &= 0, \end{aligned} \qquad [3]$$

where $\epsilon(u) = \frac{1}{2}(\nabla u + \nabla u^\top)$ is the strain rate tensor and $\nu = 10^{-3}$ is kinematic viscosity. The boundary conditions are given by $g = (1, 0)$ on the top boundary and $g = (0, 0)$ on the other three boundaries, and we seek a velocity-pressure pair $(u, p) \in V_g \times L_0^2(\Omega)$, where $V_g = \{v \in H^1(\Omega; \mathbb{R}^2) \mid v|_{\partial\Omega} = g\}$, as the pressure is only defined up to a constant. We solve Eq. (3) using a Taylor–Hood finite element discretization (12) on a $100 \times 100$ mesh, i.e., using vector-valued continuous piecewise quadratic Lagrangian elements for the velocity and piecewise linear Lagrangian elements for the pressure, implemented in the Firedrake finite element software (10) with the same solver as in Section A.

Given a known pressure term $p \in L_0^2(\Omega)$, we aim to identify the momentum and continuity equations in Eq. (3) from the observed velocity $u \in V_g$. The weak formulation of Eq. (3) becomes

$$\langle 2\nu\epsilon(u), \epsilon(\phi_u) \rangle + \langle (u \cdot \nabla)u, \phi_u \rangle = \langle p, \nabla \cdot \phi_u \rangle, \qquad [4a]$$

$$\langle \nabla \cdot u, \phi_p \rangle = 0, \qquad [4b]$$

for all test functions $(\phi_u, \phi_p) \in V_0 \times L_0^2(\Omega)$. The Dirichlet conditions are imposed strongly through the choice of the finite element nodal basis. Although the boundary function $g$ is discontinuous at the two top corners, the mismatch occurs only at isolated points (measure zero), so it does not alter the weak problem.

We first focus on the momentum equation Eq. (4a), and define the dictionary of operators $\mathcal{D} = \{\mathcal{D}_i\}_{i=1}^N$, where $N = 16$, as eight second-order spatial derivatives terms capturing viscous stress and eight nonlinear advection terms: $\{\int_\Omega \partial_j u_i \partial_k \phi_{u_i}\,\mathrm{d}x, \int_\Omega u_\ell \partial_j u_i \phi_{u_i}\,\mathrm{d}x\}$ for $i, j, k, \ell \in \{x, y\}$, with right-hand side $\mathcal{D}_0(u, p) = \int_\Omega p\nabla \cdot \phi_u\,\mathrm{d}x$. We assemble the weak forms against all nodal basis functions of the finite element discretization, leading to a system of size $80802 \times 16$, which is solved by a column-pivoted QR factorization. Here, the system is rank-deficient due to the presence of the redundant mixed derivative terms in the dictionary: $\langle \partial_x u_x \partial_y \phi_{u_x} \rangle = \langle \partial_y u_x \partial_x \phi_{u_x} \rangle$ and $\langle \partial_x u_y \partial_y \phi_{u_y} \rangle = \langle \partial_y u_y \partial_x \phi_{u_y} \rangle$, which are removed by the QR factorization. The recovered coefficients match the ground truth ones with maximum absolute error of $3.8 \times 10^{-15}$. Note that this is expected since the test functions are chosen to be the basis functions of the finite element space used to discretize Eq. (4). Hence, the discrete solution satisfies the linear system corresponding to the finite element discretization of the partial differential equation exactly, up to the tolerance of the numerical solver.

The continuity equation Eq. (4b) satisfies a homogeneous equation with zero right-hand side so its coefficients can only be recovered up to a multiplicative constant. We define the dictionary of operators $\mathcal{D} = \{\mathcal{D}_i\}_{i=1}^N$, where $N = 4$, as four first-order spatial derivatives terms capturing the divergence of the velocity field: $\{\int_\Omega \partial_j u_i \phi_p\,\mathrm{d}x\}$ for $i, j \in \{x, y\}$, with right-hand side $\mathcal{D}_0(u, p) = 0$. We assemble the weak forms using piecewise linear continuous Lagrangian elements for the test function $\phi_p$, corresponding to the dual space of our finite element discretization. This leads to a system $Gc = 0$, where $G$ has size $10201 \times 4$ and has a non-trivial null space. We aim to find a basis for the one-dimensional null space of $G$ and recover the coefficients up to a multiplicative constant by solving the following optimization problem:

$$\min_{c \in \mathbb{R}^4} \|Gc\|_2^2 \quad \text{s.t.} \quad \|c\|_2 = 1.$$

This is equivalent to finding the right singular vector of $G$ corresponding to the smallest singular value, which can be computed efficiently using the singular value decomposition (SVD). After normalizing by the coefficient of $\partial_x u_x$, we recover the continuity equation $\partial_x u_x + \partial_y u_y = 0$ with maximum absolute error of $6.7 \times 10^{-16}$.

**C. Lorenz system.** As a first ordinary differential equation example, we consider the Lorenz system (13), which is a three-dimensional system of ordinary differential equations given by

$$\dot{x} = \sigma(y - x), \qquad \dot{y} = x(\rho - z) - y, \qquad \dot{z} = xy - \beta z, \tag{5}$$

where $\sigma = 10$, $\rho = 28$, $\beta = 8/3$ so that the system features chaotic dynamics. We numerically integrate Eq. (5) from the initial condition $[1, 1, 1]^\top$ using an explicit fifth order Runge–Kutta method implemented in the `SciPy` library (14, 15), and sample the trajectory at equispaced times in the interval $[0, 40]$ with timestep $\Delta t = 10^{-3}$. We then build a dictionary consisting of all monomials up to degree two in $(x, y, z)$, along each component of the system, and a right-hand side $\mathcal{D}_0(x, y, z) = [\dot{x}, \dot{y}, \dot{z}]^\top$, so that the system decouples into three independent linear systems, one for each component of the Lorenz system. The weak formulation of Eq. (5) is obtained by multiplying the system with a test function $\phi \in \mathcal{V} = H^1([0, 40])$ and integrating by parts the right-hand side, resulting in the following weak formulation for the different components:

$$\langle c_1 + c_2 x + c_3 y + c_4 z + c_5 x^2 + c_6 xy + c_7 xz + c_8 y^2 + c_9 yz + c_{10} z^2, \phi \rangle = -\langle x, \dot{\phi} \rangle,$$
$$\langle c_1 + c_2 x + c_3 y + c_4 z + c_5 x^2 + c_6 xy + c_7 xz + c_8 y^2 + c_9 yz + c_{10} z^2, \phi \rangle = -\langle y, \dot{\phi} \rangle,$$
$$\langle c_1 + c_2 x + c_3 y + c_4 z + c_5 x^2 + c_6 xy + c_7 xz + c_8 y^2 + c_9 yz + c_{10} z^2, \phi \rangle = -\langle z, \dot{\phi} \rangle.$$

Here, we consider sinusoidal test functions (16) of the form $\phi : t \mapsto \sin(k\pi t/40)$ for $1 \leq k \leq 50$, and integrate the weak forms using the composite Simpson's rule. This leads to a system of size $50 \times 10$ for each component of the Lorenz system, which is solved by a column-pivoted QR factorization. The recovered coefficients match the ground truth ones with maximum absolute error of $5.2 \times 10^{-9}$.

**D. Duffing oscillator.** Next, we consider the forced Duffing equation modeling damped and driven oscillators, written in first-order form as

$$\dot{u} = v, \qquad \dot{v} = -\delta v - \alpha u - \beta u^3 + \gamma \cos(\omega t), \tag{6}$$

with parameters $(\alpha, \beta, \delta, \gamma, \omega) = (1, 5, 0.02, 8, 0.5)$. We integrate Eq. (6) from the initial condition $[u_0, 0]^\top$, where $u_0 \sim \mathcal{N}(0, 1)$, using an explicit fifth order Runge–Kutta method implemented in the `SciPy` library (14, 15), and store the solution $(u, v)$ at equispaced time points on the interval $[0, 100]$ with timestep $\Delta t = 10^{-3}$. We aim to recover the second order equation $\ddot{u} + \delta \dot{u} + \alpha u + \beta u^3 = \gamma \cos(\omega t)$ from the observed trajectory $(u, v)$, where $v = \dot{u}$. To do so, we select a dictionary containing all monomials up to degree three in $(u, v)$, along with the high order term $\ddot{u}$, and the right-hand side $\mathcal{D}_0 = \cos(\omega \cdot)$. We formulate the weak form of the system by multiplying the equation with sinusoidal test functions $\phi_k : t \mapsto \sin(k\pi t/100)$ for $1 \leq k \leq 50$, similarly to the Lorenz system in Section C, as

$$-\frac{1}{\gamma}\int_0^{100} v, \dot{\phi}_k\,\mathrm{d}t + \frac{1}{\gamma}\int_0^{100} \delta v + \alpha u + \beta u^3, \phi_k\,\mathrm{d}t = \int_0^{100} \cos(\omega t), \phi_k\,\mathrm{d}t$$

We discretize the integrals using Simpson's composite rule, and obtain a linear system in the dictionary coefficients. After solving the resulting linear system and normalizing by the coefficient of $\ddot{u}$, we recover the second order equation with maximum absolute error of $1.1 \times 10^{-10}$.

**E. Hankel and circulant matrices.** Finally, we consider the problem of recovering Hankel and circulant $n \times n$ matrices $A(\boldsymbol{c}^*)$ from right and left matrix-vector products: $\{(f_i, u_i = A(\boldsymbol{c}^*)f_i)\}_{i=1}^{q_R}$ and $\{(g_i, v_i = A(\boldsymbol{c}^*)^\top g_i)\}_{i=1}^{q_L}$. These matrices can be viewed as linearly parameterized operators $A(\mathbf{c}) = \sum_{i=1}^N c_i A_i$ (with $N = 2n-1$ for Hankel matrices and $N = n$ for circulant matrices), which we formulate as

$$\sum_{i=1}^N c_i \underbrace{\begin{bmatrix} A_i & & & & & \\ & \ddots & & & 0 & \\ & & A_i & & & \\ & & & A_i^\top & & \\ & 0 & & & \ddots & \\ & & & & & A_i^\top \end{bmatrix} \begin{bmatrix} f_1 \\ \vdots \\ f_{q_R} \\ g_1 \\ \vdots \\ g_{q_L} \end{bmatrix}}_{\mathcal{D}_i(u,f)} = \sum_{i=1}^N c_i \underbrace{\begin{bmatrix} A_i f_1 \\ \vdots \\ A_i f_{q_R} \\ A_i^\top g_1 \\ \vdots \\ A_i^\top g_{q_L} \end{bmatrix}}_{\mathcal{D}_i(u,f)} = \underbrace{\begin{bmatrix} u_1 \\ \vdots \\ u_{q_R} \\ v_1 \\ \vdots \\ v_{q_L} \end{bmatrix}}_{\mathcal{D}_0(u,f)}. \qquad [7]$$

Note that this formulation gives us a principled way to generalize the recovery of linearly parameterized operators to multiple evaluations of the operator on different inputs. We generate random Hankel and circulant matrices of size $n = 1000$ by sampling their coefficients from a standard normal distribution. We then sample one random vector $f_1 \sim \mathcal{N}(0, I_n)$ and compute the outputs $u_1 = A(\boldsymbol{c}^*)f_1$, and apply a column-pivoted QR linear solver to the system resulting from Eq. (7).F or circulant matrices, the resulting system is full rank with probability one for a Gaussian input vector, leading to an error of $1.6 \times 10^{-13}$ (measured in the $\ell^2$-norm) in the recovered coefficients. For Hankel matrices, one right matrix-vector product gives at most $n$ independent equations for $2n-1$ unknowns, so the system is rank deficient. Thus, we are in the case ii) of Theorem 1 in the main text, and we cannot recover the coefficients from a single input-output pair. To overcome this issue, we apply multi-shot recovery and sample two random vectors $f_1, f_2 \sim \mathcal{N}(0, I_n)$ and compute the outputs $u_1 = A(\boldsymbol{c}^*)f_1$ and $u_2 = A(\boldsymbol{c}^*)f_2$. We then assemble the system as in Eq. (7), noting that Theorem 1 applies with the operator $A(\boldsymbol{c})$ replaced by the block operator $\tilde{A}(\boldsymbol{c}) = \begin{bmatrix} A(\boldsymbol{c}) & 0 \\ 0 & A(\boldsymbol{c}) \end{bmatrix}$, and the input-output pairs $(f_1, u_1)$ and $(f_2, u_2)$, which leads to a full rank system, with coefficient recovery error of $3.5 \times 10^{-11}$. We then recovery the result reported in (17, Tab. 1), where Hankel matrices require two input-output pairs to be recovered, while circulant matrices only require one input-output pair.

## 3. Zero level set of analytic functions

For finite dimensional input spaces, one can visualize the inputs for which the parameter-to-recovery map $L_{f,u} : \mathbb{R}^N \to \mathcal{V}$ defined as $L_{f,u}(\boldsymbol{a}) = \sum_{i=1}^N a_i \mathcal{D}_i(u, f)$, for $\boldsymbol{a} \in \mathbb{R}^N$, is not injective. To do so, we select $N$ linearly independent test functionals $\psi_1, \ldots, \psi_N \in \mathcal{V}^*$ and consider the $N \times N$ matrix $G$ with entries $G_{i,j} = \psi_i(\mathcal{D}_j(u, f))$. The injectivity of $L_{f,u}$ is equivalent to the invertibility of $G$, which is characterized by the non-vanishing of the determinant $\det(G)$. Hence, the degenerate set of inputs is the zero level set of the analytic function $\phi : \mathcal{F} \to \mathbb{R}$, defined as $\phi(f) = \det(G)$, which is a thin set in the input space by Lemma 1 of the main text. In this section, we illustrate this phenomenon with the examples plotted in Figure 2 of the main text.

We begin with a two-dimensional example, where the input and residual spaces are $\mathbb{R}^2$, with dictionary functions $\mathcal{D}_1(\boldsymbol{f}) = [\cos(4f_0^2), \sin(2f_1)]^\top$ and $\mathcal{D}_2(\boldsymbol{f}) = [\sin(2f_1 + 1), \cos(2f_0)]^\top$ for $\boldsymbol{f} = [f_0, f_1]^\top \in \mathbb{R}^2$. We select canonical test functionals in $\mathbb{R}^2$ ($\boldsymbol{e}_1^\top$ and $\boldsymbol{e}_2^\top$) such that the resulting analytic function $\phi : \mathbb{R}^2 \to \mathbb{R}$ is given by

$$\phi(\boldsymbol{f}) = \det \begin{bmatrix} \cos((2f_0)^2) & \sin(2f_1 + 1) \\ \sin(2f_1) & \cos(2f_0) \end{bmatrix} = \cos((2f_0)^2)\cos(2f_0) - \sin(2f_1)\sin(2f_1 + 1), \qquad \boldsymbol{f} = (f_0, f_1) \in [-2, 2]^2.$$

In Figure 2A of the main text, we visualize the contour of $\phi$ for inputs $\boldsymbol{f} \in [-2, 2]^2$, and highlight the zero level set $\phi(\boldsymbol{f}) = 0$ in black. Lemma 1 states that the intersection of this set with a horizontal slice $f_0 = -0.3$ (plotted as dark dots) has Lebesgue measure zero in the slice, which is confirmed by the numerical experiments. Note that one could also consider the intersection with the two-dimensional plane $[-2, 2]^2$, which yields the dark zero level set curves in the figure that are of measure zero in the two-dimensional input space

Next, we consider a three-dimensional example with the dictionary functions

$$\mathcal{D}_1(\boldsymbol{f}) = [1 - 0.2f_2^2, f_1, 0]^\top, \quad \mathcal{D}_2(\boldsymbol{f}) = [2f_1, \sin(3f_0^3)\cos^2(f_2/2)(3 - f_2), 0]^\top, \quad \mathcal{D}_3(\boldsymbol{f}) = [0, 0, 1]^\top,$$

for $\boldsymbol{f} = [f_0, f_1, f_2]^\top \in \mathbb{R}^3$. Selecting canonical test functionals in $\mathbb{R}^3$, we obtain the analytic function $\phi : \mathbb{R}^3 \to \mathbb{R}$ given by

$$\phi(\boldsymbol{f}) = \det \begin{bmatrix} 1 - 0.2f_2^2 & 2f_1 & 0 \\ f_1 & \sin(3f_0^3)\cos^2(f_2/2)(3 - f_2) & 0 \\ 0 & 0 & 1 \end{bmatrix} = \left(1 - 0.2f_2^2\right)\sin(3f_0^3)\cos^2(f_2/2)(3 - f_2) - 2f_1^2,$$

where we select $\boldsymbol{f} \in [-2, 2]^3$. The zero level set of $\phi$ is visualized in Figure 2B of the main text as dark shaded surfaces, and we illustrate the intersection of this set with the plane $f_2 = 0$ that has Lebesgue measure zero in this slice.



\end{document}